\documentclass[a4paper]{amsart}

\usepackage[english]{babel}
\usepackage[utf8]{inputenc}
\usepackage{amsmath,amssymb}
\usepackage[all]{xy}
\usepackage{graphicx}
\usepackage[colorinlistoftodos]{todonotes}

\newtheorem{theorem}{Theorem}

\newtheorem{lemma}[theorem]{Lemma}

\newcommand{\Q}{{\mathbb{Q}}}
\newcommand{\Z}{{\mathbb{Z}}}
\newcommand{\R}{{\mathbb{R}}}

\title[Integral points in rational polygons]{Integral points in rational polygons: A numerical semigroup approach}

\author[G. M\'arquez--Campos]{Guadalupe M\'arquez--Campos}
\address{Departamento de \'Algebra and IMUS, Universidad de Sevilla. P.O. 1160. 41080 Sevilla, Spain.}
\email{gmarquez@us.es}
\author[J.L. Ram\'{\i}rez--Alfons\'{\i}n]{Jorge L. Ram\'{\i}rez--Alfons\'{\i}n}
\address{Institut Montpelli\'erain Alexander Grothendieck (IMAG), Université de Montpellier}
\email{jramirez@um2.fr}
\author[J.M. Tornero]{Jos\'e M. Tornero}
\address{Departamento de \'Algebra and IMUS, Universidad de Sevilla. P.O. 1160. 41080 Sevilla, Spain.}
\email{tornero@us.es}
\thanks{First and third authors were partially supported by the grant FQM--218 and P12--FQM--2696 (FEDER and FSE). The second author was supported by ANR TEOMATRO grant ANR-10-BLAN 0207 and grant ECOS-Nord M13M01.}

\date{\today}

\begin{document}

\begin{abstract}
In this paper we use an elementary approach by using numerical semigroups (specifically, those with two generators) to give a formula for the number of integral points inside a right--angled triangle with rational vertices. This is the basic case for computing the number of integral points inside a rational (not necessarily convex) polygon.
\end{abstract}

\maketitle

\section{Introduction: A little bit of history}

In recent times, impressive advances in computational combinatorics and the ever--increasing amount of applications to other branches of mathematics have made lattice--point counting a fruitful and dynamic research field. 

\vspace{.3cm}

The problem of computing the set of integral points inside plane bodies has a long history. We will see how to reduce this problem to an interated application of Theorem \ref{THR}. A milestone in this story is Pick's Theorem \cite{Pick}, from the late years of the 19th century. 

\begin{theorem}
Let $S$ be a polygon such that all of its vertices are integral, and $\mbox{int}(S)$ and $\partial S$ are, respectively, its interior and its boundary, let $A(S)$ be the area of $S$,
$$
I(S) = \# \left( \Z^2 \cap \mbox{int}(S) \right), \quad B(S) = \# \left( \Z^2 \cap \partial S \right).
$$
Then 
$$
A(S) = I(S) + \frac{B(S)}{2} -1.
$$
\end{theorem}

\vspace{.3cm}

More specifically, the question of counting lattice (in particular integral\footnote{Throughout this paper we will call a point $P \in \R^2$ {\em integral} if its coordinates lie in $\Z^2$, and similarly $P$ will be called {\em rational} if $P \in \Q^2$.}) points inside a right triangle has a long and interesting story. As early as 1922 Hardy and Littlewood \cite{HL} studied the problem of right triangles defined by the coordinate axes and a hypotenuse with irrational slope. 

\vspace{.3cm}

In the following years, the interest for the subject did not  decline. For instance, in \cite{Ehrhart} Ehrhart introduced the so--called {\em Ehrhart quasi--polynomials}, an almost ubiquous tool nowadays.  In \cite{Mordell} Mordell established a connection between lattice points problems (for the tetrahedron) and Dedekind sums while in \cite{Pommer} Pommersheim derived a more genaral formula by using techniques from algebraic geometry, see also \cite{Reeve}.  

\vspace{.3cm}

A very good and comprehensive introduction to the subject, with a good share of deep results is Beck and Robins book \cite{BR} where, in particular, the reader can find a formula, due the same authors \cite{BR1}, to compute the number of integral points inside a right triangle. Their approach is quite different from ours, they came out with a formula that uses either $n$--th primitive roots of unity or Fourier--Dedekind sums. The latter allowed them to give an efficient algorithm to calculate the number of integral points inside a right triangle since Dedekind sums with 2 variables can be computed in polynomial time. In this regard, we mention a polynomial time algorithm, due to Barnikov \cite{Barnikov1, Barnikov2} (that has significally influenced the field) that enumerates the lattice points of rational polytopes in fixed dimension.     

\vspace{.3cm}

Let us finally mention an interesting generalization of our problem. Let us call a {\em right} tetrahedron the convex set of $\R_{\geq 0}^n$ limited by the coordinate hyperplanes and a hyperplane $a_1x_1+...+a_nx_n =1$, with $a_i \in \R$. The question of counting (more precisely, bounding) the number of points in $n$--dimensional right tetrahedra has been a subject of study of S.S.T. Yau and some of his collaborators \cite{LY1,LY2,LY3,WY,XY1,XY2,XY3}, a research that produced the so--called GLY conjecture (named after its creators, A. Granville, K.P. Lin and S.S.T. Yau).

\vspace{.3cm}

This conjecture gives a lower estimate for the number of integral points in an $n$--dimensional right tetrahedron in terms of its vertex coordinates (the weak estimate) and in terms of these coordinates and the Stirling numbers (the strong estimate). The weak version was proved by Yau and Zhang \cite{YauZhang}. Asfor the strong version is concerned, the authors claim in the same paper that it has been checked computationally up to $n = 10$. The fact is that the conjecture might be checked for a particular $n$, but the state--of--the--art has not changed since. According to the authors, the case $n=10$ took weeks to be completed. This result came handy to the first and third author in \cite{MCT}. Please note that the GLY conjecture only applies to $n \geq 3$. 

\vspace{.3cm}

The main result of this paper is concerned with the following situation.

\begin{theorem}\label{THR} 
Let $a<b$ be coprime integers, $c \in \Z$. Consider the following set:
$$
T = \left\{(y_1,y_2) \in \Z_{\geq 0}^2 \; | \; a y_1 + b y_2 \leq c \right\}.
$$

Then
$$
\# T = -\frac{ab}{2} k^2 + \frac{a+b+1+2c}{2}k +  \sum_{i=0}^{\left\lfloor\frac{c - kab}{b}\right\rfloor}  \left( \left\lfloor \frac{c -kab -ib}{a} \right\rfloor +1 \right)
$$
where $k = \left\lfloor c/(ab) \right\rfloor$. Equivalently, by using the Euclidean division  $c = q\cdot ab + r$, the above formula can be expressed in the following alternative way
$$
\# T = -\frac{ab}{2} q^2 + \frac{a+b+1+2c}{2}q +  \sum_{i=0}^{\left\lfloor r/b \right\rfloor}  \left( \left\lfloor \frac{r-ib}{a} \right\rfloor +1 \right).
$$
\end{theorem}

Our technique to prove this result is based essentially in numerical semigroup invariants (the relationship between semigroups and lattice points problems has already been remarked in \cite{RA1}). Although the above formula is not polynomial in the input size (since the sum has $a$ terms in the worse case), we think that this elementary approach (avoiding Fourier analysis) might suggest a number of follow--up questions, which would be fruitful to investigate.

\vspace{.3cm}

Two such questions are shown as applications of the main result in the last section. The first one will be computing the number of integral points inside a polygon with rational vertices. The second will explore a possible application to a well--known problem in numerical semigroups: the computation of the denumerant.

\section{An interlude on numerical semigroups}

This paper relies on numerical semigroups as a fundamental tool. A numerical semigroup is a very special kind of semigroup that can be thought of as a set 
$$
\langle \, a_1,...,a_k \, \rangle = \left\{ \lambda_1 a_1 + ... + \lambda_k a_k \; | \; \lambda_i \in \Z_{\geq 0} \right\}, \mbox{ with } \gcd(a_1,...,a_k)=1.
$$

This object has been thoroughly studied in the last years, when a significant number of problems concerning the description of these semigroups and some of its more interesting invariants have been tackled. Unless otherwise stated, all proofs which are not included can be found in \cite{GS,RA}.

\vspace{.3cm}

Given a numerical semigroup $S$, there are some invariants which will be of interest for us. The most relevant will be the set of gaps, noted $G(S)$, and defined by
$$
G(S) = \Z_{\geq 0} \setminus S,
$$
which is a finite set. Its cardinality will be noted $g(S)$ and its maximum $f(S)$, the so--called Frobenius number.

\vspace{.3cm}

The Ap\'ery set of $S$ with respect to an element $s \in S$ can be defined as 
$$
Ap(S,s)= \{0, w_0, . . . ,w_{s-1}\}
$$
where $w_i$ is the smallest element in $S$ congruent with $i$ modulo $s$. 

\vspace{.3cm}

In particular, for semigroups with two generators, the invariants $g(S)$ and $f(S)$ and the relevant Ap\'ery sets are fully determined. 

\begin{lemma}\label{S2}
Let $S = \langle \, a_1, \; a_2 \, \rangle$. Then
\begin{eqnarray*}
g(S) &=& \frac{1}{2} (a_1-1)(a_2-1) \\
f(S) &=& (a_1-1)(a_2-1)-1 \\ 
Ap(S,a_i) &=& \left\{0, \, a_j, \, 2a_j, ..., \, (a_i-1)a_j \right\}
\end{eqnarray*}
\end{lemma}

\section{The number of integral points inside a right triangle}

Let us consider then a right triangle determined by the positive coordinate axes and the line 
$$
a x + b y = c, \quad a, b, c \in \Z \mbox{ and } \gcd(a,b)=1,
$$
where we will assume $a<b$, with no loss of generality. 

Take the set:
$$
T = \left\{(x,y) \in \Z_{\geq 0}^2 \; | \; ax + by \leq c \right\}
$$
and let us define the numerical semigroup associated to our triangle as $S = \langle a, \ b \rangle$. $S$ therefore verifies that its Frobenius number is $f(S) = ab -(a+b)$.

Let us perform the following partition on our set $T$:
\begin{eqnarray*}
B_0 &=& \left\{ (x, y) \in \Z_{\geq 0}^2 \; | \; ax + by \leq c, \;\; 0 \leq x < b \right\} \\
B_1 &=& \left\{ (x, y) \in \Z_{\geq 0}^2 \; | \; ax + by \leq c, \;\; b \leq x < 2b \right\} \\
\vdots && \vdots \\
B_i &=& \left\{ (x, y) \in \Z_{\geq 0}^2 \; | \; ax + by \leq c, \;\; ib \leq x < (i+1)b \right\} \\
\vdots && \vdots \\
B_k &=& \left\{ (x, y) \in \Z_{\geq 0}^2 \; | \; ax + by \leq c, \;\; kb \leq x \right\}
\end{eqnarray*}
where $k= \lfloor c/(ab) \rfloor$. 

\vspace{.5cm}
\begin{figure}
\centering
\includegraphics[scale=0.6]{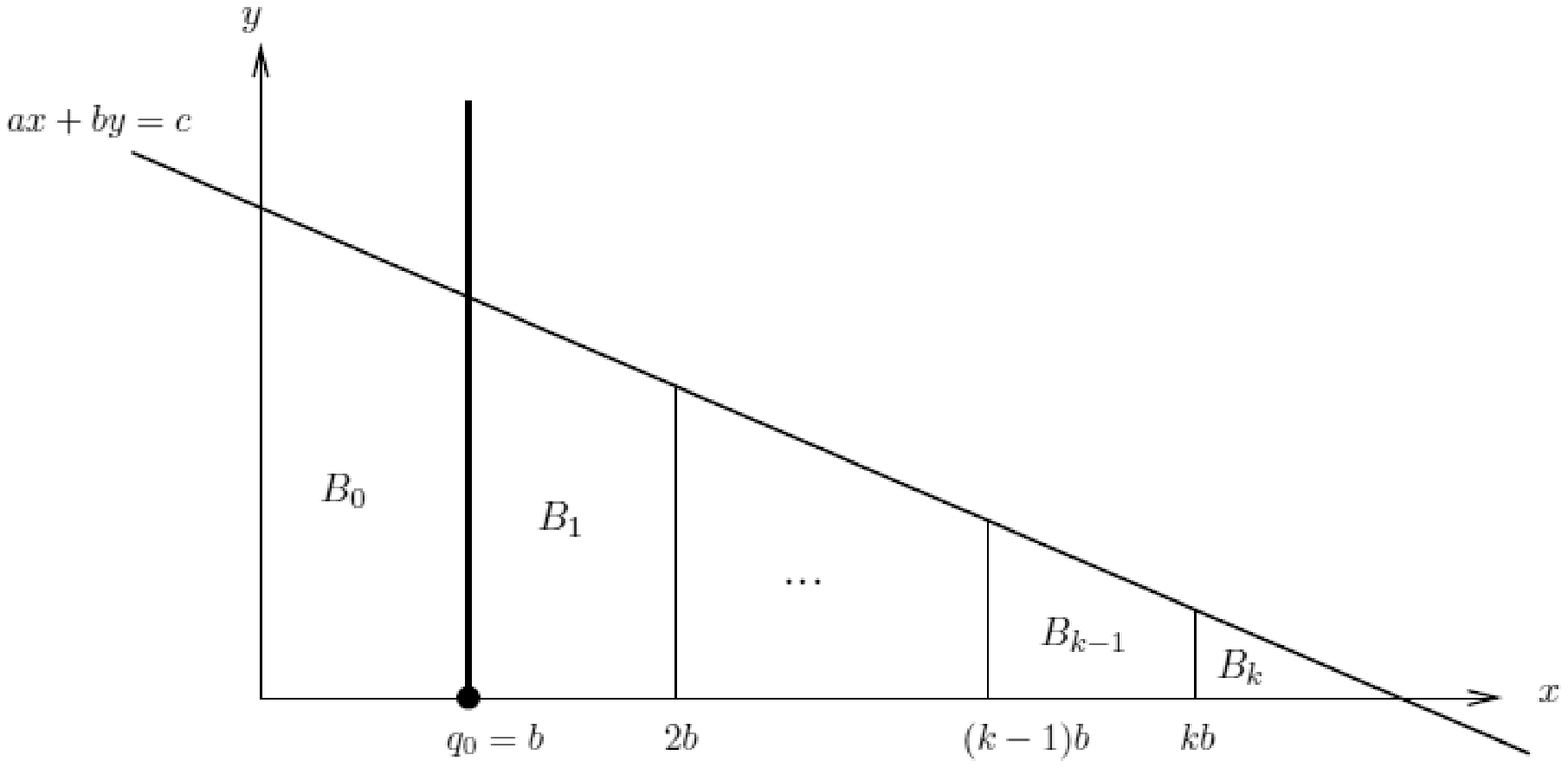}
\end{figure}
\vspace{.5cm}

As our aim is to find $\# T$, and it is plain that:
$$
\# T = \# B_0 + \# B_1 + ... + \# B_{k-1} + \# B_k,
$$
we can reduce our problem to that of finding $\# B_i$, for $i=0,...,k$.

\begin{lemma}\label{B0}
Under the previous assumptions, if $k>0$,
$$
\# B_0 = \frac{a+b-ab+1}{2} + c = \frac{1-f(S)}{2} + c
$$
\end{lemma}

\begin{proof}
We will actually show that
$$
S \cap [0,c] \stackrel{1:1}{\longleftrightarrow} B_0.
$$

Given a pair $(x,y) \in B_0$ we define its corresponding element in $S \cap [0,c]$ to be $n=ax+by$. 

In the same way, given $n \in S \cap [0,c]$ it is clear that we must have a representation $n=ax+by$ and we can in fact assume $0 \leq x < b$ (if otherwise, we can move part of $ax$ into the $by$ summand until $x<b$). 

Let us assume that we have such a representation (that is, with $0 \leq x < b$) and we will prove that then the pair $(x,y)$ must be unique, which will establish the bijection. Should we have
$$
n = ax_0+by_0 = ax_1+by_1, \quad \mbox{ with } 0<x_0,x_1<b,
$$
we must have 
$$
a(x_0-x_1) = b(y_1-y_0)
$$
and, as $\gcd(a,b)=1$, this means $b|(x_1-x_0)$, which in turn implies $x_0=x_1$. Note that the case $x_0=x_1=0$ leads directly to $y_0=y_1$ as desired.

Note that $k>0$ is equivalent to $c \geq ab$, which also yields $c > f(S)$. Therefore, after Lemma \ref{S2},
\begin{eqnarray*}
\# B_0 = \displaystyle \# \big( S \cap [0,c] \big) &=& \# \Big( S \cap [0, f(S)] \Big) + c - f(S) \\
&=& \frac{ab -(a+b)+1}{2} + c - \big( ab-(a+b) \big) \\
&=& \frac{a+b-ab+1}{2} + c.
\end{eqnarray*}
\end{proof}

Simple as it is, this case is the basic argument for the whole process. Now, if we want to compute $\# B_1$, we just move our triangle, so that $(b,0)$ is now at the origin. Similarly, the line $ax+by=c$ is moved, as in the picture:

\vspace{.5cm}

\includegraphics[scale=0.6]{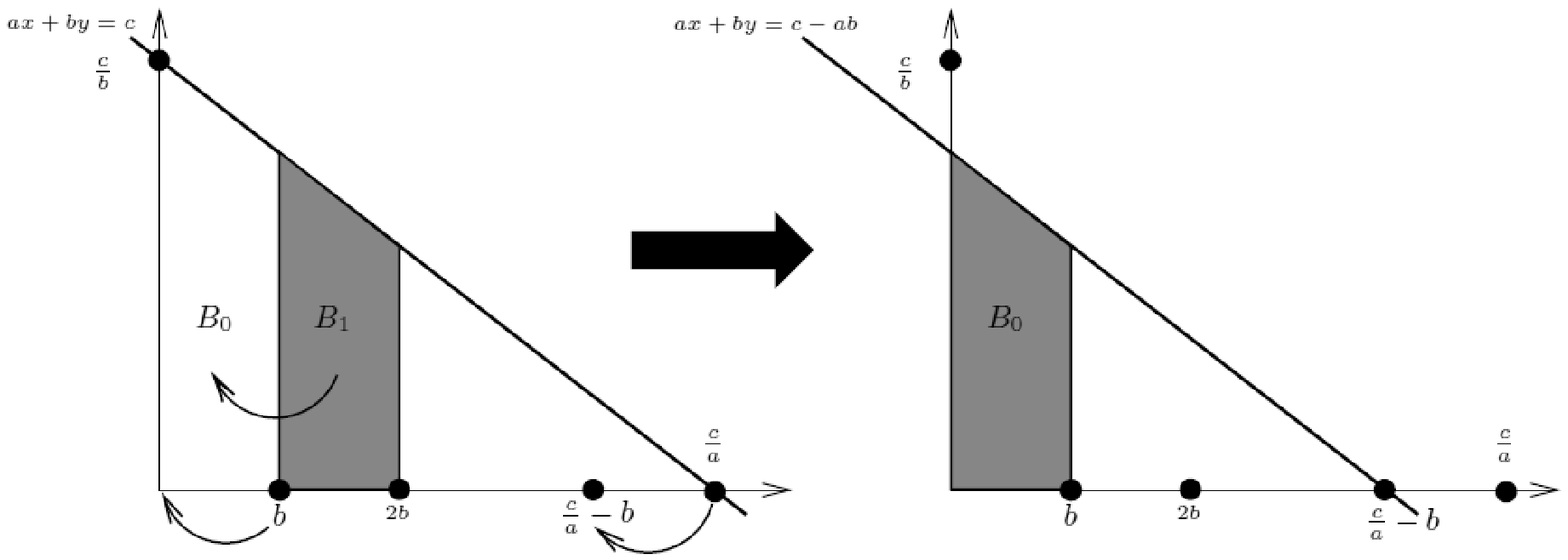}

\vspace{.5cm}

So, with a little abuse of notation, let us redefine:
$$
B_1 = \left\{ (x, y) \in \Z_{\geq 0}^2 \; | \; ax + by \leq c_1, \;\; 0 \leq x < b \right\}
$$
where $c_1 = c-ab$ (obviously this only makes sense if $c > ab$). Assuming $k>1$ we have, following the same way:
\begin{eqnarray*}
\# B_1 = \#  \left( S \cap [0, c_1] \right) &=& \# \Big( S \cap [0, c - ab] \Big) \\
&=& \frac{ab -(a+b)+1}{2} + c -ab - (ab-(a+b) + 1) \\
&=& \frac{a+b-3ab+1}{2} + c.
\end{eqnarray*}

\vspace{.3cm}

We can of course go along the same lines for computing $\# B_i$ for $i= 1,...,k-1$, where $k= \lfloor c/(ab) \rfloor$, rewriting $c_i = c -i ab$, whenever $c > iab$ and we will find:
\begin{eqnarray*}
\# B_i = \# \left( S \cap [0, c_i] \right) &=& \# \Big( S \cap [0, c - iab] \Big) \\
&=& \frac{ab -(a+b)+1}{2} + c - iab - (ab-(a+b)) \\
&=& \frac{(a+b)-(1+2i)ab + 1}{2} + c
\end{eqnarray*}

\vspace{.3cm}

We have then arrived at the nutshell of the problem: the set $B_k$. After we have moved it to the origin, we have our renamed $B_k$:
$$
B_k = \left\{ (x,y) \in \Z_{\geq 0}^2 \; | \; ax + by \leq c_k \right\}.
$$

Now we might have $c_k < ab-(a+b)$. So we cannot proceed in the same way as before. We do know $c_k = c - kab$, that is, $c_k = c \mod ab$, and from Lemma \ref{S2} we also know:
$$
Ap(S,a) = \{ 0, b, 2b, ..., (a-1)b \}
$$
and therefore 
$$
\{ w \in Ap(S,a) \; | \; w \leq c_k \} = \left\{ib \; \Big| \;  i= 0,1,...,\left\lfloor \frac{c_k}{b} \right\rfloor \right\}
$$

On the other hand, if $i \in \{0,...,\lfloor c_k/b \rfloor \}$, we have
$$
ib+ja \leq c_k \; \Longleftrightarrow \; j \leq \left\lfloor \frac{c_k -ib}{a} \right\rfloor,
$$
and then
\begin{eqnarray*}
S \cap [0, c_k] &=& \{ ib +ja \leq c_k \; | \; i,j \in \Z_{\geq 0} \} \\
&=& \left\{ ib +ja \; \Big| \;  i \in \left\{0,..., \left\lfloor \frac{c_k}{b} \right\rfloor \right\}, \;\; j \leq \left\lfloor \frac{c_k - ib}{a} \right\rfloor \right\} \\
&=&  \sum_{i=0}^{\left\lfloor c_k/b \right\rfloor} \left( \left\lfloor \frac{c_k -ib}{a} \right\rfloor +1 \right)
\end{eqnarray*}

Adding up all of these computations, we arrive to our result. In the previous conditions:
\begin{eqnarray*}
\# T &=& \# B_0 + \# B_1 + ... + \# B_i + ... + \# B_{k-1} + \# B_k \\
&=& \left(\frac{a+b-ab+1}{2} + c\right) + ... + \left(\frac{(a+b)-(1+2i)ab +1}{2} + c\right) + \\
&& \quad \quad \quad +...+  
\sum_{i=0}^{\left\lfloor c_k/b \right\rfloor} \left( \left\lfloor \frac{c_k -ib}{a} \right\rfloor +1 \right) \\
&=& \sum_{i=0}^{k-1}  \left(\frac{(a+b)-(1+2i)ab + 1}{2} + c\right) + \sum_{i=0}^{\left\lfloor c_k/b \right\rfloor} \left( \left\lfloor \frac{c_k -ib}{a} \right\rfloor +1 \right) \\
&=& -\frac{ab}{2} k^2 + \frac{a+b+1+2c}{2}k +  \sum_{i=0}^{\left\lfloor\frac{c - kab}{b}\right\rfloor}  \left( \left\lfloor \frac{c -kab -ib}{a} \right\rfloor +1 \right)
\end{eqnarray*}
where $k = \left\lfloor c/(ab) \right\rfloor$. 
Or equivalently, by using the Euclidean division: $c = q\cdot ab + r$ we obtain
$$
\# T = -\frac{ab}{2} q^2 + \frac{a+b+1+2c}{2}q +  \sum_{i=0}^{\left\lfloor r/b \right\rfloor}  \left( \left\lfloor \frac{r-ib}{a} \right\rfloor +1 \right).
$$
This proves Theorem \ref{THR}.

\vspace{.3cm}

\noindent {\bf Example. } 
Let us do a simple example to illustrate the process, considering the triangle defined by the line $3x + 7y = 46$.

Following Theorem \ref{THR} we have $k = \left\lfloor \frac{46}{3\cdot 7} \right\rfloor =2$, that is,
$$
\# \ T =- \frac{3\cdot7}{2} 2^2 + \frac{3+7+1+2\cdot 46}{2}\cdot 2 +  \sum_{i=0}^{\left\lfloor\frac{46 - 2\cdot 3\cdot 7}{7}\right\rfloor}  \left( \left\lfloor \frac{46 -2 \cdot 3\cdot 7 -i \cdot 7}{3} \right\rfloor +1 \right)= 63.
$$

Let us see the actual counting:

\vspace{.5cm}

\includegraphics[scale=0.66]{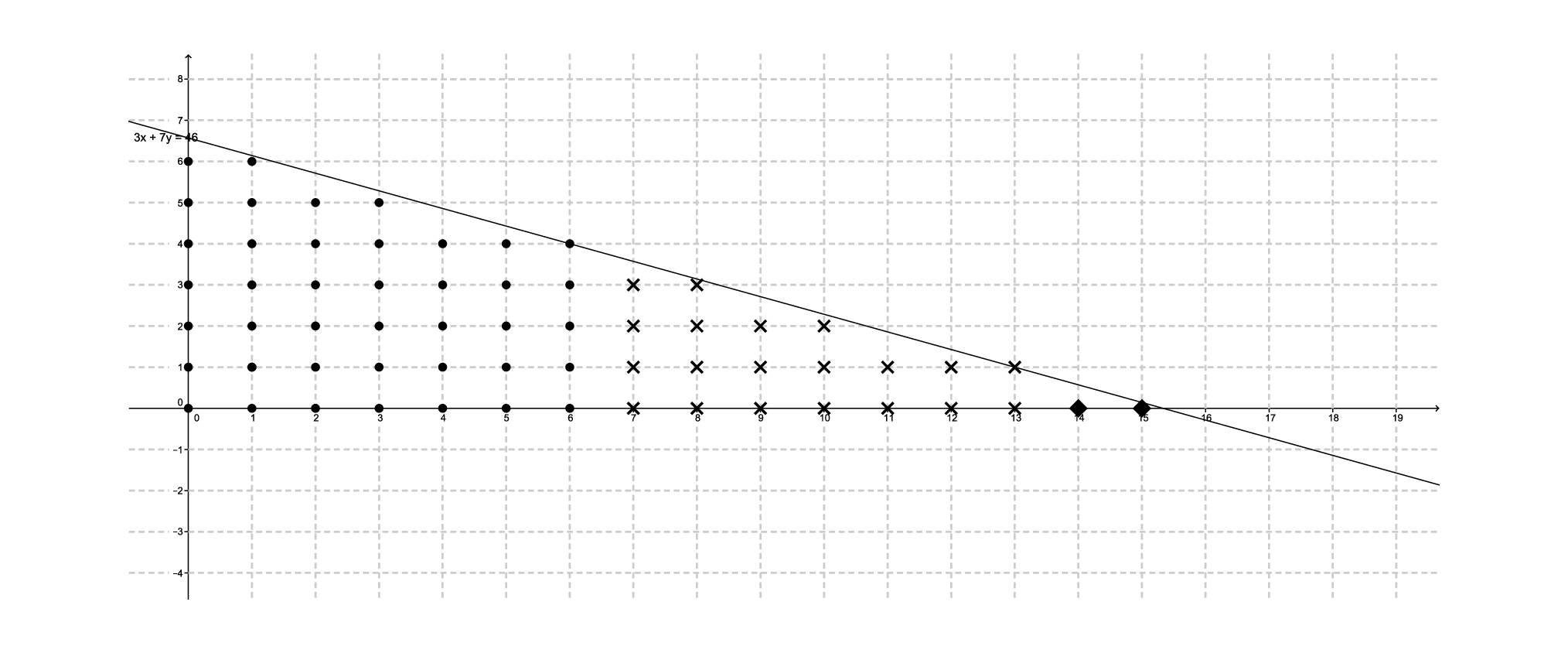}

\vspace{.5cm}

In the picture we have put different symbols for the different sets $B_j$, following the process. Round points correpond to $B_0$. There are $41$ of them, as predicted by the formula:
$$
\# B_0 = \frac{a+b-ab+1}{2}+c = \frac{3+7-3\cdot 7+1}{2} + 46 = 41.
$$

Crossed points correspond to points in $B_1$: 
$$
\# B_1 = \frac{a+b-3ab+1}{2}+c = \frac{3+7-3\cdot 3\cdot 7+1}{2} + 46 = 20.
$$

And finally the square points are those of $B_2$:
$$
\# B_2 = \sum_{i=0}^{\left\lfloor\frac{c_k}{b}\right\rfloor} \left( \left\lfloor \frac{c_k -ib}{a} \right\rfloor +1 \right) = \left\lfloor \frac{46-42}{3} \right\rfloor +1 = 2.
$$

\section{Applications (I): Integral points in general triangles}

Our aim in this section is to give a result (in some sense in the spirit of Pick's theorem) to compute the number of integral points inside a generic triangle. This, of course, is enough to compute the number of integral points inside polygons (not necessarily convex) defined by rational vertices. But we will not address here the problem of dividing a polygon into triangles \cite{BKOS,KKT}.

The important point here is that, in order to compute the number of integral points of such a polygon, it suffices with rectangles and right triangles of a particular type. Of course rectangles can be divided in two right triangles, but this brings no substantial computational simplification to our problem.

\vspace{.3cm}

From now on, a rectangle whose sides are parallel to the coordinate axes will be called a {\em stable} rectangle. Computing the number of integral points inside a stable rectangle is very easy.

\begin{lemma}
Let $\alpha_1<\beta_1$ and $\alpha_2<\beta_2$ be real numbers. Let $R \subset \R^2$ be the stable rectangle with vertices $(\alpha_1,\alpha_2)$ and $(\beta_1,\beta_2)$. Then
$$
\# (R \cap \Z^2) =  \left( \lfloor \beta_1 \rfloor - \lceil \alpha_1 \rceil + 1 \right) \cdot  \left( \lfloor \beta_2 \rfloor - \lceil \alpha_2 \rceil + 1 \right)  
$$
\end{lemma}

Similarly, a right triangle whose orthogonal sides are parallel to the coordinate axes will be called a {\em stable} right triangle. The computation of the number of integral points inside a stable right triangle can be easily achieved from Theorem \ref{THR}, as we shall show now.

Assume we have such a triangle $T$ defined by rational vertices ($A$ is the vertex at which $T$ has its right angle):
$$
A = (\alpha,\beta), \quad B = (\alpha,\gamma), \quad C = (\delta,\beta),
$$
and we can assume, up to symmetry, $\alpha<\gamma$, $\beta<\delta$. Furthermore, it is clear that, as for counting integral points is concerned, we can substitute
$$
\alpha \longmapsto \lceil \alpha \rceil, \quad \quad \beta \longmapsto \lceil \beta \rceil,
$$
and the number of integral points does not change by traslations of integral vectors, hence we can in fact assume $A = (0,0)$.

\vspace{.5cm}

\begin{center}
\includegraphics[scale=0.8]{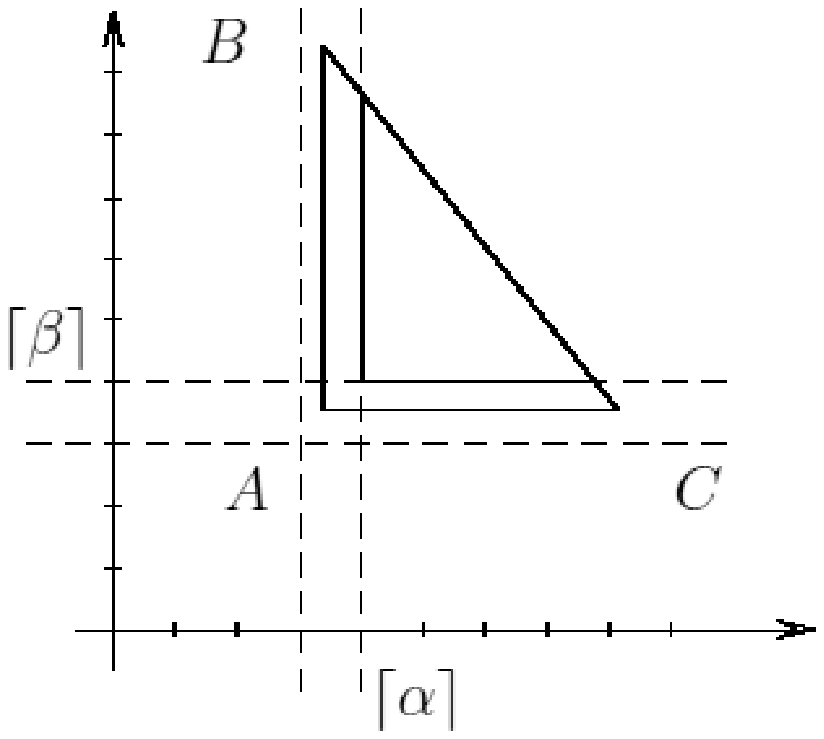}
\end{center}

\vspace{.5cm}

After translating the triangle $T$, we can write the hypothenuse of the new triangle as
$$
ax+by = c,
$$
where we can assume $a,b,c \in \Z$, $\gcd(a,b,c)=1$. That is, $B = (0,c/b)$, $C = (c/a,0)$. Mind that $a$ and $b$ need not to be coprime in this set--up. But, in order to apply Theorem \ref{THR}, we should have $\gcd(a,b)=1$.

\vspace{.3cm}

Let $d= \gcd(a,b)$. There are two possibilities:

\begin{enumerate}
\item $d=1$. We can apply Theorem \ref{THR}.
\item $d>1$. As $\gcd(d,c)=1$ clearly there are no integral points in the hypothenuse, as any such point $(x,y)$ must verify $ax+by \in \Z d$. So we are going to change the triangle $T$, moving the hypothenuse in a parallel way towards $A$, until the new triangle has at least some chance to have one integral point in the hypothenuse. This is the triangle $T'$ defined by
$$
ax+by=c',
$$
where $c'= \lfloor c/d \rfloor d$. It is clear that with such construction $\# (T \cap \Z^2) = \# (T' \cap \Z^2)$ and so 
$\# (T \cap \Z^2)$ can be computed counting the number of integral points in the stable triangle defined by
$$
\frac{a}{d}x + \frac{b}{d} y = \frac{c'}{d} = \left\lfloor \frac{c}{d} \right\rfloor,
$$
that is, using Theorem \ref{THR}.
\end{enumerate}

\vspace{.3cm}

As we shall see later, in related problems (as, for instance, the problem of computing the number of integral points inside a polygon), we might have to compute the number of points in the sides of a right stable triangle. To do this, note that:

\begin{itemize}
\item The number of integral points in the orthogonal sides can be easily computed: in the segment limited by $(a,b)$ and $(a,c)$, say with $b<c$ (any other case is obviously symmetric), the number of integral points is
$$
 \left( \lfloor c \rfloor - \lceil b \rceil + 1 \right) \cdot 1_\Z (a),
$$
where $1_\Z (a)$ is the indicator function of $\Z$. 

\item The number of integral points in the line $ax+by=c$ is the number of representations of $c$ inside the semigroup $S = \langle a, \ b \rangle$. This number is known as the denumerant of $c$ in $S$ \cite{RA}. 

It is easy to see that this denumerant has to be either $\lfloor c/(ab)\rfloor$ or $\lfloor c/(ab)\rfloor +1$, that is, $k$ or $k+1$ with the notation of the previous section, assuming $c \in S$ (obviously it is $0$ otherwise). This is because $c$ (or $c_i$) must be representable and from Lemma \ref{B0} there must exist exactly one representation whose coefficients are in $B_0$ (respectively $B_i$). This holds true for $B_0$,...,$B_{k-1}$ but not necessarily for $B_k$ (because $c_k$ might not be in $S$), hence the two possible cases. More precisely, we have the following result (see \cite{Popo} for the original proof in Romanian, \cite{BCS} for a shorter and easier one):
\end{itemize}

\begin{theorem}
Under the previous asumptions, let $a'$ and $b'$ be the only integers veryfing
$$
0 < a'< b, \quad a \cdot a' = -c \mod b
$$
$$
0 < b'< a, \quad b \cdot b' = -c \mod a
$$

Then the denumerant of $c$ in $S$ is given by
$$
d(c;a,b) = \frac{c+a\cdot a' + b \cdot b'}{ab} -1
$$
\end{theorem}

In the sequel, when we are interested in the points inside a right triangle $T$ we will write $\#T^{int}$ and if we want the points excluding only the hypothenuse we will denote this cardinal by $\# T^{hyp}$.

\vspace{.3cm}

Finally, then, let us show how to compute the integral points in a generic triangle using stable rectangles and right triangles (and hence, Theorem \ref{THR}). Assume we are given a triangle $T$, given by
$$
A_1 = (\alpha_1,\beta_1), \quad A_2 = (\alpha_2,\beta_2), \quad A_3 = (\alpha_3,\beta_3),
$$
and let us call
$$
x_0 = \min_{i=1,2,3} \alpha_i, \quad x_1 = \max_{i=1,2,3} \alpha_i;  
$$
$$
y_0 = \min_{i=1,2,3} \beta_i, \quad y_1 = \max_{i=1,2,3} \beta_i;
$$

We consider the (stable right) rectangle $R_T$ determined by 
$$
\Big\{ \, V_1 = (x_0,y_0), \; V_2 = (x_1,y_0), \; V_3 = (x_0,y_1), \; V_4 = (x_1,y_1) \, \Big\}.
$$

It is clear that it must hold $V_i = A_j$ for some $i=1,2,3,4$ and $j=1,2,3$. Depending on the relative position of the maximal vertices, the situation must be one of these four (up to rotation and flip, if necessary):

\begin{itemize}
\item Only one of the vertices of $R_T$ is a vertex of $T$. Then it is clear that:
$$
\# T = \# R_T - \# T_1^{hyp} - \# T_2^{hyp} - \# T_3^{hyp}.
$$

\vspace{.5cm}

\includegraphics[scale=0.4]{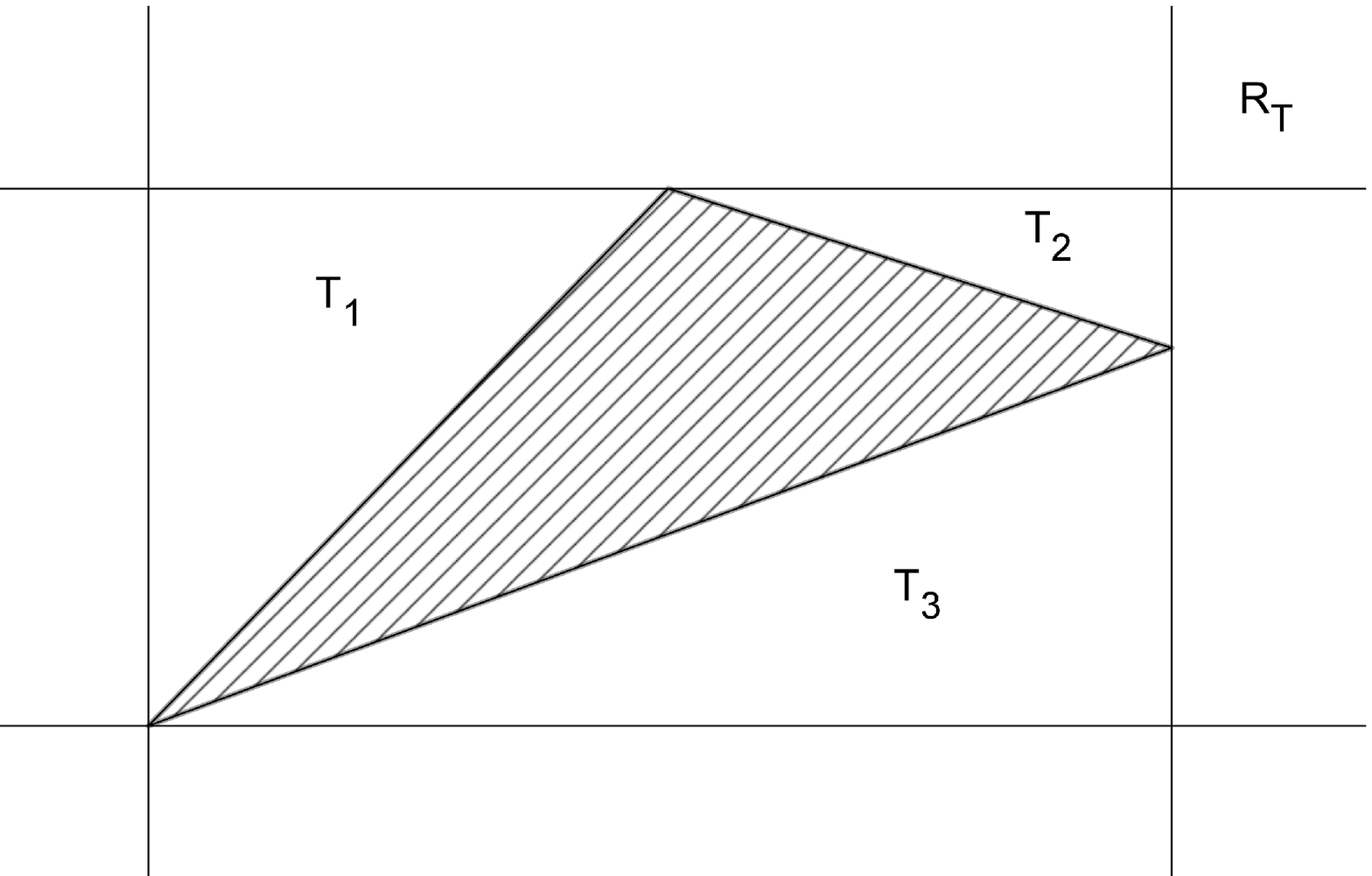}

\vspace{.5cm}

\item Two adjacent vertices of $R_T$ are vertices of $T$. Then:
$$
\# T = \# R_T - \# T_1^{hyp} - \# T_2^{hyp}.
$$

\vspace{.5cm}

\includegraphics[scale=0.4]{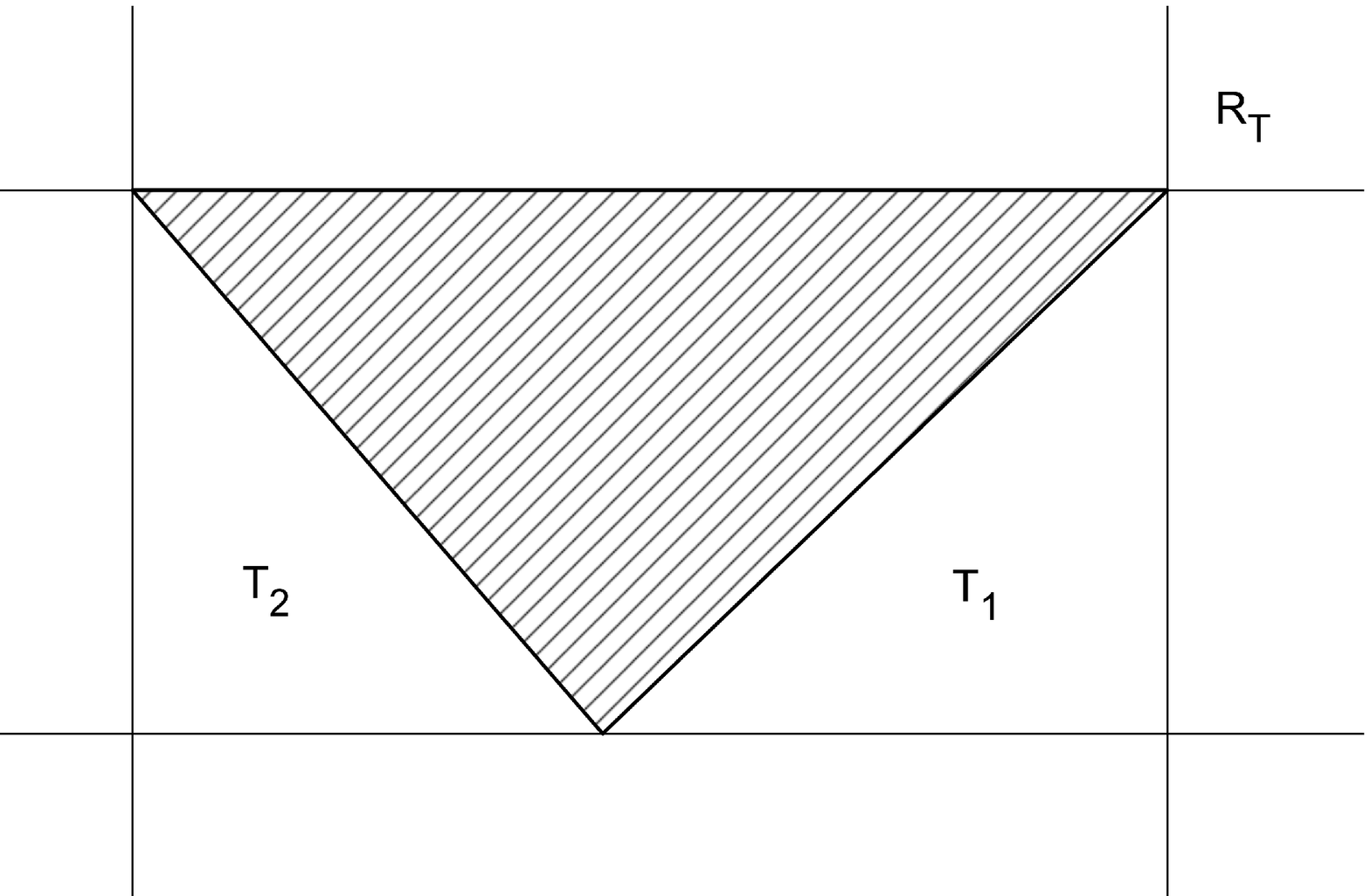}

\vspace{.5cm}

\item Two non--adjacent vertices of $R_T$ are vertices of $T$. Then (one of many possibilities):
$$
\# T = \# R_T - \# T_1^{hyp}  - \# T_2^{hyp} - \# T_3^{hyp} - \# T_4 - \# T_5^{int} + 1_{\Z^2}(P) + 1_{\Z^2}(Q).
$$

\vspace{.5cm}

\includegraphics[scale=0.4]{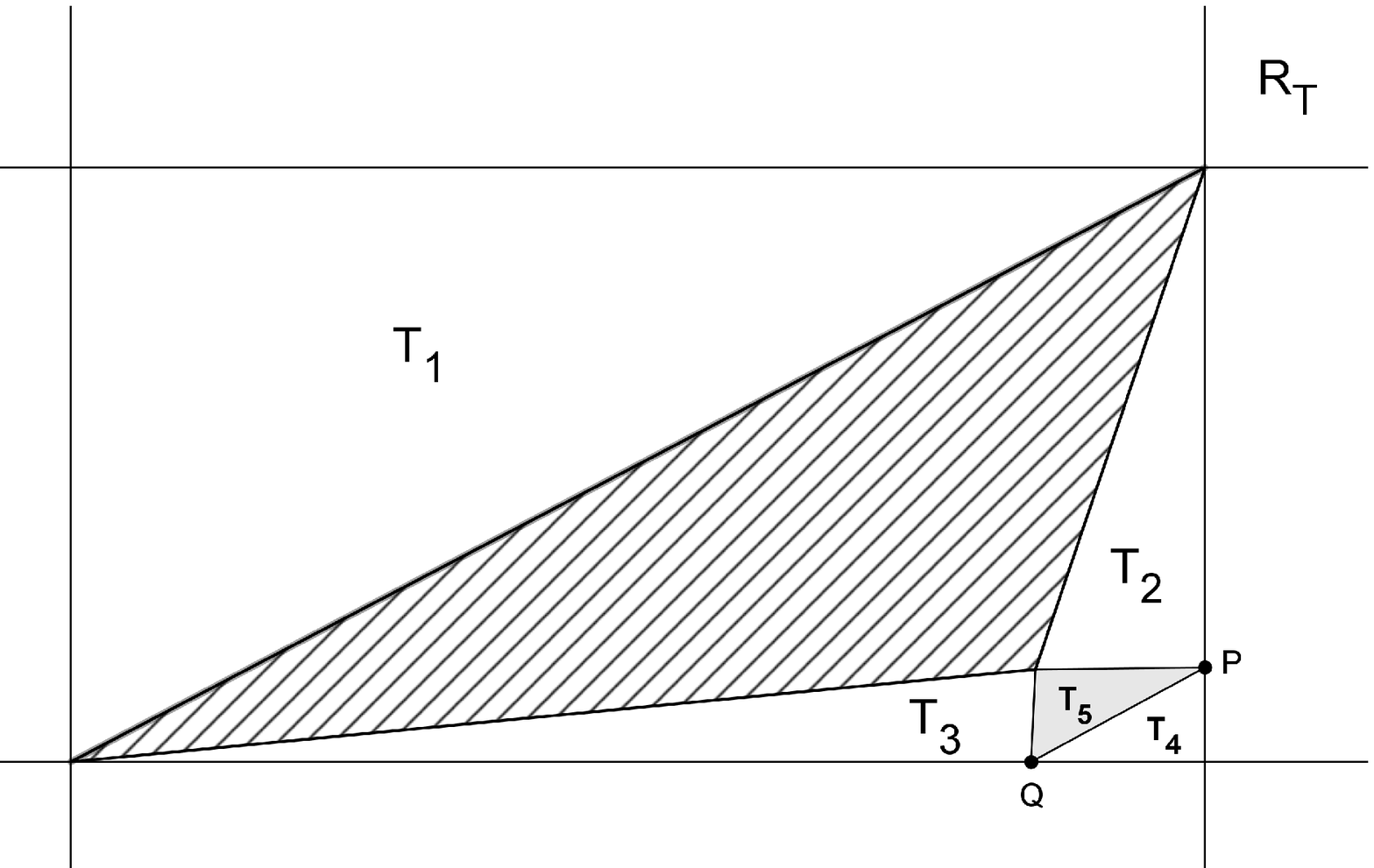}

\vspace{.5cm}

\item The three vertices of $T$ are vertices of $R_T$. This is the trivial case, as $T$ is a right stable triangle.
\end{itemize}

\section{Applications (II): On the denumerant of a numerical semigroup with $3$ generators}

The result on rational polygons and its reduction to rectangles and right triangles could be generalized to an $n$--dimensional set--up. However, the actual formulas are not easy enough so as to give a tight result. We give a first idea on how this could be done by computing the number of points in a stable right tetrahedron (the definition is the obvious one).

\begin{theorem}
Let $T (a_1,a_2,a_3,b) \subset \R^3$ be the tetrahedron defined by
$$
T (a_1,a_2,a_3,b) = \left\{ (x_1,x_2,x_3) \; | \; x_i \geq 0, \; a_1x_1+a_2x_2+a_3x_3 \leq b \right\}
$$
where we are assuming $a_1 < a_2 < a_3$, $\gcd(a_1,a_2)=1$.

For $i=0,...,\lfloor b/a_3 \rfloor$ define $q_i$ and $r_i$ by the Euclidean division:
$$
b-a_3i = q_i(a_1a_2) + r_i.
$$

Then 
\begin{eqnarray*}
\# \left(T (a_1,a_2,a_3,b) \cap \Z^3 \right) &=& \sum_{i=0}^{\lfloor b/a_3 \rfloor} \Bigg( -\frac{ab}{2} q_i^2 + \frac{a+b+1+2(b-a_3i)}{2}q_i + \\
&& \quad\quad + \sum_{j=0}^{\lfloor r_i/b \rfloor}  \left( \left\lfloor \frac{r_i-ja_2}{a_1} \right\rfloor +1 \right) \; \Bigg).
\end{eqnarray*}
\end{theorem}

\begin{proof}
The formula is just the result of adding the number of points in every right triangle $T(a_1,a_2,a_3,b) \cap \{ x_3 = i \}$ for $i = 0,...,\lfloor b/a_3 \rfloor$.
\end{proof}

\vspace{.5cm}
\begin{center}
\includegraphics[scale=0.5]{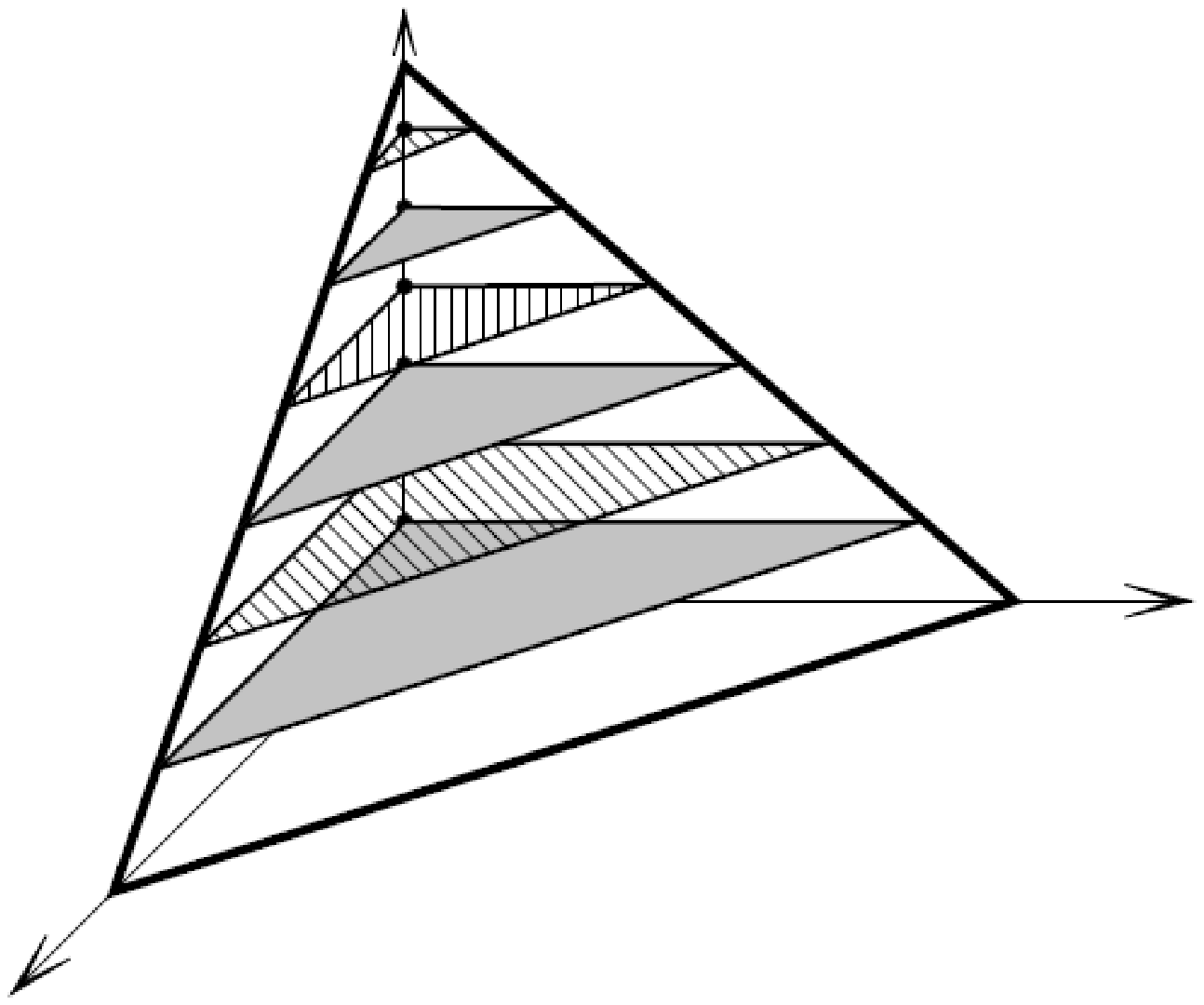}
\end{center}
\vspace{.5cm}

The condition $\gcd(a_1,a_2)=1$ can obviously be substituted by $\gcd(a_1,a_3)=1$ or $\gcd(a_2,a_3)=1$ if necessary. If none of this conditions is met, like in the tetrahedron defined by
$$
6x_1+10x_2+15x_3=21,
$$
for instance, then some of the right triangles have to be {\em adjusted} as we did in the previous section. This is not a difficulty when programming, so to say, but the general formula gets a lot more complicated. We have tried to get a compact version of this, but this effort has been unsuccesful so far.

This result can be handy when trying to compute the denumerant function we introduced above.

Easy as it is to define, the denumerant is a very elusive function which has proved elusive to compute even in cases with $3$ generators (see \cite[Chapter 4]{RA}). With the previous result one can give a formula, not very sophisticated though. Simply note that 
\begin{eqnarray*}
d \left( a;a_1,a_2,a_3 \right) &=& \# \Big\{ (x_1,x_2,x_3) \in \Z_{\geq 0} \; | \; a_1x_1+a_2x_2+a_3x_3 = a \Big\} \\
&=& \# \left(T (a_1,a_2,a_3,a) \cap \Z^3 \right) - \# \left(T (a_1,a_2,a_3,a-1) \cap \Z^3 \right).
\end{eqnarray*}

And then, from the previous result, one can obtain the desired formula.

\vspace{.5cm}
\begin{center}
\includegraphics[scale=0.5]{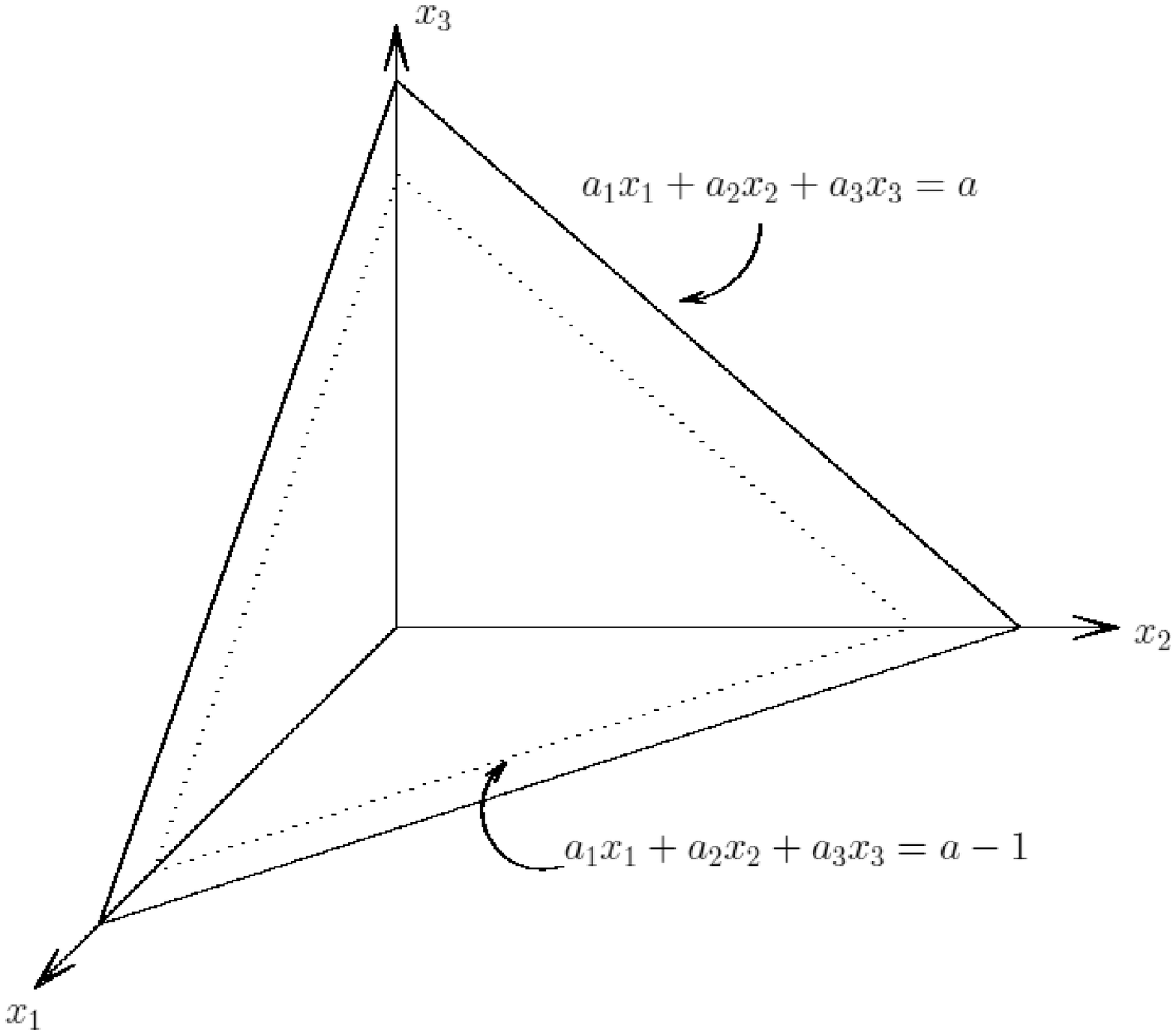}
\end{center}
\vspace{.5cm}

\section{Acknowledgements}

The authors warmly acknowledge the help and advice of J.C. Rosales. 

\vspace{.3cm}

The final stage of this work was completed during a stay of the first author in the Institut Montpelli\'erain Alexander Grothendieck, Universit\'e de Montpellier, during the spring term of 2013, thanks to the grant P08--FQM--3894 (Junta de Andaluc\'{\i}a).

\end{document}